# An elementary proof of Fermat's Last Theorem

YURI ARENBERG

**Abstract**

We show that an elementary proof of Fermat's Last Theorem (FLT) exists. Our paper also extends the scope of FLT from integers to all rational numbers. Starting with the $a^n + b^n = c^n$ equation, where $a, b$ and $c$ are real nonzero numbers and $n$ is an integer equal to or greater than 2, its $a, b$ and $c$ components are converted to expressions which are linear in the same term. This procedure also yields an $n - 1$ degree equation. To prove FLT, it is necessary to demonstrate that, for $n \geq 3$, this equation has no solutions in rational numbers. Our analysis reveals that, if some solution in rational numbers were to exist, it would be accompanied by another solution in rational numbers, a counterpart to the original. It is then shown that, when $n \geq 3$, the structure of the $n - 1$ degree equation does not sustain the coexistence of the counterpart solutions in rational numbers. FLT is the consequence of this result.

**Mathematics Subject Classification (2010)**    11D41    11Z05    12E1

## 1  Introduction

Fermat's Last Theorem (FLT) asserts that it is impossible to separate a power greater than the second into two like powers. Circa 1637, Fermat noted in his copy of Arithmetica that he proved this proposition. However, he apparently never publicized his claim of a general proof; and only an indication of his proof of the theorem's fourth power case survived.

In 1994, Andrew Wiles (collaborating with Richard Taylor) succeeded in proving FLT. Rooted in the Taniyama-Shimura conjecture (now the modularity theorem), his ultra-modern proof could not have been conceived several decades earlier, much less during Fermat's time. Thus, with the issue of whether it is even possible to prove or disprove FLT settled, the question of the existence a 17[th] century proof remained unresolved. The purpose of our paper is to fill this gap.

In what follows, we show that an elementary proof of FLT exists. Our paper also extends the scope of this theorem from integers to all rational numbers. Starting with the $a^n + b^n = c^n$ equation, where $a, b$ and $c$ are real nonzero numbers and $n$ is an integer equal to or greater than 2, its $a, b$ and $c$ components are converted to expressions which are linear in the same term. This procedure also yields an $n - 1$ degree equation. To prove FLT, it is necessary to demonstrate that, for $n \geq 3$, this equation has no solutions in rational numbers. Our analysis reveals that, if some solution in rational numbers were to exist, it would be accompanied by another solution in rational numbers, a counterpart to the original; and we uncover a rule governing the relationship between these solutions. It is then shown that, when $n \geq 3$, the structure of the $n - 1$ degree equation does not support the coexistence of the counterpart solutions in rational numbers. FLT is the consequence of this result.

The paper is organized as follows. Section 2 sets forth the framework of this study. Section 3 is devoted to examining the $n = 2$ and $n = 3$ cases and identifying the reason why $a^3 + b^3 = c^3$ cannot hold in rational numbers. A general proof of FLT is presented in Section 4.

## 2 The framework

Consider a triple $\{a, b, c\}$ of real nonzero numbers satisfying the following equation:

$$a^n + b^n = c^n, \tag{2.1}$$

where $n$ is an integer and $n \geq 2$. Without any loss of generality, it will be assumed that $a > 0$, $b > 0$, $c > 0$, $c > a$, $c > b$.

Let $c = a + d$ and $a + b = c + f$, where $d$ and $f$ are real numbers, $d > 0$, $f > 0$, $a > f$ and $b > f$. Then $b = d + f$. Substituting $d + f$ for $b$ and $a + d$ for $c$ in (2.1), expanding $(a + d)^n$ and $(d + f)^n$ and rearranging yields

$$f^n + nf^{n-1}d + \tfrac{1}{2}n(n-1)f^{n-2}d^2 + \cdots + \tfrac{1}{2}n(n-1)f^2 d^{n-2} + nfd^{n-1}$$

$$-na^{n-1}d - \tfrac{1}{2}n(n-1)a^{n-2}d^2 - \cdots - \tfrac{1}{2}n(n-1)a^2 d^{n-2} - nad^{n-1} = 0. \tag{2.2}$$

Let $f = gd$ and $a = ghd$, where $g$ and $h$ are real numbers, $g > 0$ and $h > 1$. Substituting $ghd$ for $a$ and $gd$ for $f$ in (2.2), dividing each component of the resulting expression by $gd^n$ and collecting the terms leads to the following equation:

$$g^{n-1} - n(h^{n-1} - 1)g^{n-2} - \tfrac{1}{2}n(n-1)(h^{n-2} - 1)g^{n-3} - \cdots$$

$$- \tfrac{1}{2}n(n-1)(h^2 - 1)g - n(h-1) = 0. \tag{2.3}$$

It must be satisfied simultaneously with

$$a = ghd, \ b = (g+1)d, \ c = (gh+1)d. \tag{2.4}$$

This transformation introduces the possibility that (2.3) may have roots which are negative numbers. Without any loss of generality, let us primarily focus on solutions in positive numbers.

Expression (2.4) reveals that (2.1) can hold in rational numbers only if all members of a triple $\{d, g, h\}$ on which it is based are rational numbers. Because $a$, $b$ and $c$ are linear in $d$, one can always make $d$ a rational number. Hence, to prove FLT, it is sufficient to show that (2.3) cannot be satisfied if $\{g, h\}$ is a pair of rational numbers and $n \geq 3$.

**Proposition 2.1.** For each solution to (2.3) in real numbers, there must be another solution in real numbers, a counterpart to the original.

*Proof.* Instead of letting $c = a + d$, the equivalents of (2.3) and (2.4) could have obtained by setting $c = b + d$. Consequently, if (2.3) holds for a pair $\{g_1, h_1\}$ of some triple $\{d_1, g_1, h_1\}$ of real numbers there must also be a counterpart triple $\{d_2, g_2, h_2\}$ of real numbers whose members $g_2$ and $h_2$ also satisfy (2.3) and, additionally, in view of (2.4), fulfill these three requirements:

$$g_1 h_1 d_1 = (g_2 + 1)d_2, \quad (g_1 + 1)d_1 = g_2 h_2 d_2, \quad (g_1 h_1 + 1)d_1 = (g_2 h_2 + 1)d_2. \tag{2.5}$$

From (2.5), the following results can be deduced:

$$g_1 = \frac{1}{h_2 - 1}, \quad g_2 = \frac{1}{h_1 - 1}. \tag{2.6}$$

Clearly, if a pair $\{g_1, h_1\}$ were to exist in rational numbers, $\{g_2, h_2\}$ would also be a pair of rational numbers; and the solution to (2.3) based on one of these pairs would be transformable to that based on its counterpart pair via (2.6). Our strategy in proving FLT is to show that, for any $n \geq 3$, the underlying structural relationship between $h$ and $g$ in (2.3) precludes the coexistence of the counterpart solutions in rational numbers. This conclusion is reached by demonstrating that, if $n \geq 3$, the transformation from one solution to (2.3) to its counterpart solution via (2.6) cannot even be accomplished in real numbers.

### 3 The n=2 versus the n=3 case

Before proceeding to a general proof, it would be instructional to consider the case of $n = 2$ and that of $n = 3$ and pinpoint the reason why solutions in rational numbers are possible in the former, but not in the latter one.

If $n = 2$, (2.3) reduces to

$$g = 2(h - 1). \tag{3.1}$$

The existence of solutions in rational numbers is not in question in this case. However, to illustrate our strategy, we shall offer a formal proof.

*Proof.* Let $\{g_1, h_1\}$ and $\{g_2, h_2\}$ be some two counterpart pairs of positive real numbers, each of them satisfying (3.1). Then $g_1 = 2(h_1 - 1)$ and $g_2 = 2(h_2 - 1)$. Now, once, say, $h_1$ has been specified and $g_1$ ascertained from (3.1), $h_2$ and $g_2$ are determinable from (2.6). To prove that solutions in rational numbers exist, it is necessary to show first that $g_1 = 2(h_1 - 1)$ is transformable to $g_2 = 2(h_2 - 1)$ via (2.6). The former expression may be written as

$$\frac{1}{h_2 - 1} = 2(h_1 - 1). \tag{3.2}$$

Now, to obtain $g_2$, the left hand side (LHS) term of (3.2) must be multiplied by $\frac{h_2 - 1}{h_1 - 1}$. Then, if the product of multiplication of the right hand side (RHS) term of (3.2) and this ratio were to match the RHS term of $g_2 = 2(h_2 - 1)$, the transformation would be accomplished. There is indeed a match.

To complete the proof, one must show that (3.2) can hold in rational numbers. This requirement is satisfied by making $h_1$ a rational number. It will be shown that, for any $n \geq 3$, the structure of (2.3) does not accommodate transformations in real numbers via (2.6). This finding at once rules out the existence solutions in rational numbers when $n \geq 3$.

In the $n = 3$ case, (2.3) becomes

$$g^2 = 3(h^2 - 1)g + 3(h - 1). \tag{3.3}$$

**Theorem 3.1** (The $n = 3$ Case of Generalized FLT for Rational Numbers). *There does not exist a triple of positive rational numbers $\{a *, b *, c *\}$ satisfying the equation $a^3 + b^3 = c^3$*. The first published proof for positive integers is due to Euler. For an excellent description and analysis of his proof and other major developments on the FLT front prior to the 20th century, see [1].

*Proof.* Assume that, contrary to Theorem 3.1, there exists some triple $\{a *, b *, c *\}$ of positive rational numbers for which $a^3 + b^3 = c^3$ holds. Hence there must also be two counterpart pairs of positive rational numbers, $\{g_1 *, h_1 *\}$ and $\{g_2 *, h_2 *\}$, each of them satisfying (3.3). A necessary condition for their existence mandates that there must be two counterpart pairs of positive real numbers, $\{g_1, h_1\}$ and $\{g_2, h_2\}$, whose members, in view of (3.3), are related as

$$g_1^2 = 3(h_1^2 - 1)g_1 + 3(h_1 - 1) \tag{3.4}$$

and

$$g_2^2 = 3(h_2^2 - 1)g_2 + 3(h_2 - 1), \tag{3.5}$$

and that (3.4) be transformable to (3.5) via (2.6). Before attempting to perform this transformation, let us consider a special case.

**Lemma 3.1.** If $g_1 = \frac{1}{h_1 - 1}$ (which, in view of (2.6), implies that $h_1 = h_2$ and $g_1 = g_2$), solutions in rational numbers are not possible.

*Proof.* To reach this conclusion, substitute $\frac{1}{h_1 - 1}$ for $g_1$ in (3.4) and, after letting $h_1 = \frac{k}{l}$, where $k$ and $l$ ($k > l > 0$) are relatively prime integers, obtain $l^3 = 6k(k - l)^2$. This expression cannot hold in integers. It can be verified that, if $g_1 = \frac{1}{h_1 - 1}$, solutions in rational numbers do not exist in any $n \geq 2$ case.

In view of (2.6), one may write (3.4) and (3.5) as

$$\frac{1}{(h_2 - 1)^2} = \frac{3(h_1^2 - 1)}{h_2 - 1} + 3(h_1 - 1) \tag{3.6}$$

and

$$\frac{1}{(h_1 - 1)^2} = \frac{3(h_2^2 - 1)}{h_1 - 1} + 3(h_2 - 1). \tag{3.7}$$

Now, once $h_1$ has been specified, the structure of (3.6) becomes embedded in that of (3.7). Consequently, a valid (3.6) to (3.7) transformation must preserve this structure while producing a matching result.

Applying the multiplicative operator $\left(\frac{h_2 - 1}{h_1 - 1}\right)^2$ to every term on each side of (3.6) yields

$$\frac{1}{(h_1 - 1)^2} = \frac{3(h_1 + 1)(h_2 - 1)}{h_1 - 1} + \frac{3(h_2 - 1)^2}{h_1 - 1}. \tag{3.8}$$

This operator does not only secure a match of the (3.7) and (3.8) LHS terms, but also ensures that the structure of (3.6) is preserved in (3.8). However, it does not bring about the required match of the (3.7) and (3.8) RHS terms. In search for this match, one may be tempted to perform additive operations on the two RHS terms of (3.8). Indeed a matching result would have been seemingly obtained if $\frac{3(h_2-1)(h_2-h_1)}{h_1-1}$ were added to the first term and subtracted from the second. The validity of such manipulations is addressed by the lemma below.

**Lemma 3.2.** Additive operations in connection with transformations of expressions based some polynomial equation of degree 2 or higher are structure preserving only if these expressions are evaluated at this equation's roots. The proof is straightforward; and we shall omit it.

Now, $\frac{1}{h_1-1}$ is not a root of (3.6). In addition to $\frac{1}{h_2-1}$, this equation's other root is $-3(h_1-1)(h_2-1)$. Therefore, all additive manipulations of the RHS terms of (3.8) must be ruled out.

It follows that, for the (3.6) to (3.7) transformation to be possible, the RHS terms of (3.8) must at once, in their current form, match those of (3.7). Because this requirement has not been satisfied, it is necessary to conclude that $\{g_1 *, h_1 *\}$ and $\{g_2 *, h_2 *\}$ do not exist. This contradictory result negates the assumption of the existence of $\{a *, b *, c *\}$, completing the proof of Theorem 3.1. As seen, the lack of solutions in rational numbers in this case is ascribable to the fact that the structure of (3.3) does not even allow transformations in real numbers via (2.6).

## 4 The general case

**Theorem 4.1** (Generalized FLT for Rational Numbers). *There does not exist a triple of positive rational numbers $\{a *, b *, c *\}$ satisfying the equation $a^n + b^n = c^n$ for $n \geq 3$. See [2, 3] for details of Wiles' proof.*

*Proof.* The proof follows the reasoning espoused in Section 3. Assume that, contrary to Theorem 4.1, there exist some triple $\{a *, b *, c *\}$ of positive rational numbers for which (2.1) holds and hence two counterpart pairs of positive rational numbers, $\{g_1 *, h_1 *\}$ and $\{g_2 *, h_2 *\}$, each of them satisfying (2.3). To confirm that a necessary condition for the existence of $\{g_1 *, h_1 *\}$ and $\{g_2 *, h_2 *\}$ has been satisfied, one must show that, for some two counterpart pairs of positive real numbers, $\{g_1, h_1\}$ and $\{g_2, h_2\}$, whose members, in view of (2.3), are related as

$$g_1^{n-1} = n(h_1^{n-1} - 1)g_1^{n-2} + \frac{1}{2}n(n-1)(h_1^{n-2} - 1)g_1^{n-3} + \cdots$$
$$+ \frac{1}{2}n(n-1)(h_1^2 - 1)g_1 + n(h_1 - 1) \tag{4.1}$$

and

$$g_2^{n-1} = n(h_2^{n-1} - 1)g_2^{n-2} + \frac{1}{2}n(n-1)(h_2^{n-2} - 1)g_2^{n-3} + \cdots$$
$$+ \frac{1}{2}n(n-1)(h_2^2 - 1)g_2 + n(h_2 - 1), \tag{4.2}$$

(4.1) can be transformed to (4.2) via (2.6).

In view of (2.6), (4.1) and (4.2) become

$$\frac{1}{(h_2-1)^{n-1}} = \frac{n(h_1^{n-1}-1)}{(h_2-1)^{n-2}} + \frac{n(n-1)(h_1^{n-2}-1)}{2(h_2-1)^{n-3}} + \dots$$

$$+ \frac{n(n-1)(h_1^2-1)}{2(h_2-1)} + n(h_1-1) \tag{4.3}$$

and

$$\frac{1}{(h_1-1)^{n-1}} = \frac{n(h_2^{n-1}-1)}{(h_1-1)^{n-2}} + \frac{n(n-1)(h_2^{n-2}-1)}{2(h_1-1)^{n-3}} + \dots$$

$$+ \frac{n(n-1)(h_2^2-1)}{2(h_1-1)} + n(h_2-1). \tag{4.4}$$

The first step in the (4.3) to (4.4) transformation process entails multiplying every term on each side of the former by $\left(\frac{h_2-1}{h_1-1}\right)^{n-1}$. It yields

$$\frac{1}{(h_1-1)^{n-1}} = \frac{n(h_1^{n-2}+h_1^{n-3}+\dots+h_1+1)(h_2-1)}{(h_1-1)^{n-2}}$$

$$+ \frac{n(n-1)(h_1^{n-3}+h_1^{n-4}+\dots+h_1+1)(h_2-1)^2}{2(h_1-1)^{n-2}} + \dots$$

$$+ \frac{n(n-1)(h_1+1)(h_2-1)^{n-2}}{2(h_1-1)^{n-2}} + \frac{n(h_2-1)^{n-1}}{(h_1-1)^{n-2}} \tag{4.5}$$

and produces a match of the LHS terms of (4.4) and (4.5). However, these expressions' RHS terms are yet to be matched. Now, consistent with our arguments in the previous section, on the one hand, solutions in rational numbers are not possible if $g_1 = \frac{1}{h_1-1}$. On the other hand, unless $\frac{1}{h_1-1}$ is one of the roots of (4.3), the structure of this equation and hence that of (4.4) would be altered if additive manipulations were performed on the RHS terms of (4.5). Consequently, for the (4.3) to (4.4) transformation to be valid, the RHS terms of (4.5) must at once, in their current form, match those of (4.4).

A careful inspection of (4.4) and (4.5) reveals that the required match occurs only when $n = 2$. This finding forces a conclusion that, if $n \geq 3$, $\{g_1 *, h_1 *\}$ and $\{g_2 *, h_2 *\}$ do not exist, invalidating the existence of $\{a *, b *, c *\}$ assumption for any $n \geq 3$ case and completing the proof of Theorem 4.1. The culprit for the absence of solutions in rational numbers in the general case is the structure of (2.3). Unless $n = 2$, it does not even accommodate transformations in real numbers via (2.6).

Department of Economics, Brooklyn College and The Graduate Center, City University of New York (CUNY)

E-mail address: YuriA@brooklyn.cuny.edu

**Acknowledgements and Dedication**

The author would like to thank Kenneth Ribet, Barry Mazur and Brian Conrad for their advice on procedural matters and Marc Fox, Gary Testa and, especially, Paul Goldberg for their encouragement. Above all, the author is grateful to Peter D. Loeb and Jacob Sturm.  Their insightful questions and comments regarding this paper's earlier draft led to a significant improvement of the quality of exposition. Professor Sturm's further suggestions shaped this version the paper. Any of its shortcomings is attributable solely to the author.  In memory of my mother.